\documentclass[12pt]{amsart}
\usepackage[margin=1in]{geometry}
\newcommand{\Z}{{\mathbb Z}}

\newcommand{\F}{{\mathbb F}}
\newcommand{\ord}{\mbox{ord}}
\newcommand{\pp}{{\mathfrak p}}
\newcommand{\calF}{{\mathcal F}}
\newcommand{\calX}{{\mathcal X}}
\newcommand{\calB}{{\mathcal B}}
\newcommand{\calC}{{\mathcal C}}
\newcommand{\calS}{{\mathcal S}}
\newcommand{\qq}{{\mathfrak q}}
\newcommand{\Gal}{\mbox{Gal}}
\newcommand{\id}{\mbox{id}}
\newcommand{\calE}{{\mathcal E}}
\font\te=eufm10
\newcommand{\mte}[1]{\mbox{\te {#1}}}
\newcommand{\isom}{\simeq}
\newtheorem{theorem}{Theorem}[section]
\newtheorem{lemma}[theorem]{Lemma}
\newtheorem{corollary}[theorem]{Corollary}
\newtheorem{proposition}[theorem]{Proposition}

\theoremstyle{definition}
\newtheorem{definition}[theorem]{Definition}
\newtheorem{question}[theorem]{Question}

\theoremstyle{remark}
\newtheorem{remark}[theorem]{Remark}

\newtheorem{notation}[theorem]{Notation}

\begin{document}%
\bibliographystyle{alpha}%
\title{Undecidability in function fields of positive
  characteristic}%
\author{Kirsten Eisentr\"ager}
\address{Department of Mathematics\\
The Pennsylvania State University\\
University Park, PA 16802, USA.}
\author{Alexandra Shlapentokh}
\address{Department of Mathematics\\
East Carolina University\\
Greenville, NC 27858, USA.}
\begin{abstract}
  We prove that the first-order theory of any function field $K$ of
  characteristic $p> 2$ is undecidable in the language of rings
  without parameters. When $K$ is a function field in one variable
  whose constant field is algebraic over a finite field, we can also
  prove undecidability in characteristic 2. The proof uses a result
  by Moret-Bailly about ranks of elliptic curves over function fields.
\end{abstract} 


\thanks{K.\ Eisentr\"ager was partially supported by a
  National Science Foundation postdoctoral fellowship. A.~Shlapentokh
  was partially supported by NSF grants DMS-0354907 and DMS-0650927.
  }
\maketitle

\section{Introduction}\label{intro}
The current investigation started as an attempt by the authors to
resolve Hilbert's Tenth Problem for all function fields of positive
characteristic.  Hilbert's Tenth Problem in its original form was to
find an algorithm to decide, given a polynomial equation
$f(x_1,\dots,x_n)=~0$ with coefficients in the ring $\Z$ of integers,
whether it has a solution with $x_1,\dots,x_n \in \Z$.
Matiyasevich (\cite{Mat70}), building on earlier work by Davis,
Putnam, and Robinson (\cite{DPR61}), proved that no such algorithm
exists, i.e.\ Hilbert's Tenth Problem is undecidable.

Since then, analogues of this problem have been studied by asking the
same question for polynomial equations with coefficients and solutions
in other recursive commutative rings. Perhaps the most important
unsolved question in this area is Hilbert's Tenth Problem over the
field of rational numbers.

 The function field analogue turned out to
be much more tractable.
We know that Hilbert's Tenth Problem for the function field $k$ of a
curve over a finite field is undecidable.  This was proved by Pheidas
for $k=\F_q(t)$ with $q$ odd~(\cite{Ph3}), and then extended to all
global function fields in \cite{V,Sh13,Eis}.
We also have undecidability of Hilbert's Tenth Problem for certain
function fields over possibly infinite constant fields of positive
characteristic~(\cite{Sh15, Sh30,Eis,K-R1}).  The results of
\cite{Eis} and \cite{Sh15} also generalize to higher transcendence
degree (see \cite{Sh24}) and give undecidability of Hilbert's Tenth
Problem for finite extensions of $\F_q(t_1,\dots,t_n)$ with $n\geq 2$.
In \cite{Eispositive} the problem was shown to be undecidable for
finite extensions of $k(t_1, \dots, t_n)$ with $n \geq 2$ and $k$
algebraically closed of odd characteristic.

So all known undecidability results for Hilbert's Tenth Problem in
positive characteristic either require that the constant field not be
algebraically closed or that we are dealing with a function field in
at least 2 variables.
The big open question that remains is whether Hilbert's Tenth Problem
for a one-variable function field over an algebraically closed field
of constants is undecidable.

The current methods for proving undecidability of Hilbert's Tenth
Problem for function fields $K$ of positive characteristic $p$ usually
require showing that the following sets are existentially definable in
the language of rings : $\{(x,x^{p^s}):x \in K, s \in \Z_{\geq 0}\}$
and $\{x\in K: \ord_{\pp}x \geq~0\} \mbox{ for some nontrivial prime
}\pp \mbox{ of } K$.
In this paper we show that we can existentially define one of these
sets for a large class of fields: we will prove that the set of $p$-th
powers is existentially definable in {\em any} function field $K$ of
characteristic $p>2$ whose constant field has transcendence degree at
least one over $\F_p$.

By a {\em function field (in $n$ variables)} over a field $F$ we mean a field
$K$ containing $F$ and $n$ elements $x_1, \dots, x_n$, algebraically
independent over $F$, such that $K/F(x_1,\dots, x_n)$ is a finite
algebraic extension. The algebraic closure of $F$ in $K$ is called the
constant field of $K$, and it is a finite extension of $F$.

Given the present difficulties of showing that Hilbert's Tenth
Problem, or, equivalently, the existential theory of an arbitrary
function field of positive characteristic is undecidable, one can also
consider a weaker result, namely proving the undecidability of the
first-order theory of these fields.
Duret showed that the first-order theory of function fields (in $n$
variables) over algebraically closed fields of positive characteristic
is undecidable (\cite{Duret2}).  In \cite{Cherlin} Cherlin showed that
the first-order theory of $F(t)$ is undecidable for infinite perfect
fields $F$ of positive characteristic. In \cite{Ph7} Pheidas extended
this result to rational function fields $F(t)$ for any field $F$ of
characteristic $p \geq 5$, but he had to add a transcendental parameter
$t$ to the ring language to prove undecidability.

In this paper we generalize Duret's and Pheidas' results and prove
that the first-order theory in the language of rings without
parameters of {\it any} function field over a field of characteristic
greater than 2 is undecidable. In the case the field of constants is
algebraic over a finite field, we can also treat the case of
characteristic 2.

The paper is organized as follows. We first show that the
first-order theory for function fields $K$ of positive characteristic
is undecidable in the language of rings with finitely many parameters.
This is done by defining a model of the nonnegative integers with
addition and multiplication in $K$.  Then, using a result of R.\
Robinson, we show that this also gives us the undecidability of the theory
of $K$ in the language of rings without parameters.  The details of
this argument are discussed in Section~\ref{coeff}.

We endeavored to make the presentation as uniform as possible across
all the different types of fields. Thus, in the second section of the
paper we first show that in order to establish the first-order
undecidability of a function field of characteristic $p>0$, it is
enough to show that $p$-th powers of a specific field element are
first-order definable.  To define $p$-th powers of a specific element,
the techniques we use depend on the constant field. When the constant
field is algebraic over a finite field we generalize equations that
have previously been used in \cite{Sh15,Eis} to reduce the problem to
the rational function field case. This is done in Section~\ref{td1}.
When the constant field has transcendence degree $\geq 1$ over $\F_p$,
things are more complicated. In Section~\ref{highertd} we show that
the $p$-th powers we want to define occur as $x$-coordinates of the
$K$-rational points on a certain elliptic curve. We then use a theorem by
Moret-Bailly about the rank of elliptic curves in extensions function
fields to reduce to the rational function field case. The theorem by
Moret-Bailly was also used in \cite{Eis2,Eispadic,MB3} to obtain
undecidability results. Moret-Bailly's theorem only holds in odd
characteristic, and so for higher transcendence degree we obtain
undecidability for function fields of odd positive characteristic.

\pagebreak

\section{Using $p$-th powers to construct a  model of the positive
  integers}
\subsection{Statement of results}
The main result that we will prove in the first four sections is the
following:
\begin{theorem}\label{maintheorem}
  Let $K$ be a function field of characteristic $p>2$ or a function
  field in one variable of characteristic 2 whose constant field is
  algebraic over $\F_2$. There exists a finite set of parameters
  $\{z_1,\dots z_k\} \subseteq K$ (depending on $K$) such that first-order
  theory of $K$ is undecidable in the language of rings augmented by
  $\{z_1,\dots,z_k\}$.
\end{theorem}
From the result in Section~\ref{coeff} we obtain a strengthening of
Theorem~\ref{maintheorem}:
\begin{theorem}\label{maintheorem2}
  Let $K$ be as in Theorem~\ref{maintheorem}. Then the first-order
  theory of $K$ is undecidable in the language of rings without
  parameters.
\end{theorem}
\subsection{Idea of proof}
In this section we show that to prove Theorem~\ref{maintheorem}, it
is enough to define $p$-th powers of an element in the field with at
least one simple pole or zero. To prove that this is enough we use a
result of J.~Robinson that shows how to define multiplication of
positive integers in terms of addition and divisibility.

\begin{remark}
  Since we are only interested in the function field $K$ and not the
  underlying field $F$, we can always replace $F$ with $F(x_2, \dots,
  x_n)$ and view a function field $K/F$ in $n$ variables as a function
  field $K/F(x_2, \dots, x_n)$ in one variable.  So in the following
  all the function fields we consider will be function fields in one
  variable.
\end{remark}
\begin{notation}%
\label{not:K}
Let $K$ be a function field (in one variable) of positive
characteristic $p$ over a field of constants $F$. Let $F_0$ be the
algebraic closure of a finite field in $F$. Let $t \in K \setminus F$,
and let $n=[K:F(t)]$.
\end{notation}%

The following result is due to J.~Robinson (see \cite{Rob1}).
\begin{lemma}
  There exists a first-order formula $\calF$ in the language
  $\langle\Z_{>0},+,\mid \rangle$ such that for integers $k,m,n$, we have $k=mn
  \iff \calF(k,m,n)$. Here $a \mid b$ means ``$a$ divides $b$'' for
  positive integers $a,b$.
\end{lemma}

An immediate corollary of this lemma is the fact that the first-order
theory of $\langle \Z_{>0},+,\mid \rangle$ is undecidable.
So to prove the undecidability of the first-order theory of $K$ it is
enough to construct a model of the positive integers with addition and
divisibility in $K$.

We say that we have a {\em model} of
$\langle\Z_{>0},+,\mid \rangle$ in $K$ if there is a bijection $\phi:
\Z_{>0} \to D$ between $\Z_{>0}$ and a definable subset $D$ of $K^d$
(for some $d \geq 1$), such that the graphs of $+$ and $\mid$ on $D$
induced by $\phi$ correspond to definable subsets of $D^3$ and $D^2$,
respectively.

As we will see in Theorem~\ref{thm:model} below, we can construct a
model of $\langle\Z_{>0},+,\mid\rangle$ if we can define $p$-th powers
of a specific element.
\subsection{From $p$-th powers of a special element to arbitrary
  $p$-th powers}
The known definitions of $p$-th powers in general are produced in the
following manner: first define $p$-th powers of a specific element,
then use $p$-th powers of this element to produce $p$-th powers of
arbitrary elements. 

We have to distinguish two cases: the case where the constant field is
perfect and when it is not perfect.  The function fields $K$
considered in Section~\ref{td1} have constant fields which are
algebraic over finite fields. These constant fields are always
perfect. However, this will not be necessarily true of the constant
fields in Section~\ref{highertd}.  The constant field will be relevant
in two places: in the definition of the $p$-th powers of the special
element in Lemma~\ref{le:powersoft2} and also in the proof of
Proposition~\ref{prop:needt} below.

\begin{proposition}%
\label{prop:needt}
Let $K$ be a function field of positive characteristic.
Suppose the set $p(K,t):=\{x \in K: \exists s \in \Z_{> 0}, x=
t^{p^s}\}$ is definable in $K$ for some $t \in K \setminus F$ which is
not a $p$-th power in $K$. Then the following subset of $K^2$ is also
definable in $K$:
\[%
\calX(K):=\{(x, x^{p^s}): s \in \Z_{> 0}, x \in K\}.
\]%
\end{proposition}%

\begin{proof}
  When the constant field is perfect we can follow the arguments of
  Section 8.4 of~\cite{Sh34} which covers all positive
  characteristics.  While the function fields considered in
  \cite{Sh34} are global function fields, the only condition that is
  really required is that the constant field is perfect.  The main
  ingredient in the proof is that we can define a global derivation
  and use it to determine whether certain elements of the field have
  simple zeroes and poles.

  When the field of constants is not perfect we obtain the result from Lemma
  2.10 and Corollary 2.11 of \cite{Sh15}.  
In both cases we might need
  to enlarge the field of constants to avoid ramifying valuations as
  zeros and poles of elements whose $p$-th powers we define.  These
  constant field extensions are covered by Proposition
  \ref{le:convert} and Lemma \ref{le:ramify} from the appendix. So the
  proposition holds for any field of constants.
\end{proof}

Both \cite{Sh34} and \cite{Sh15} also prove the following corollary
which will be needed below.
\begin{corollary}%
\label{cor:samepower}%
 Let $t \in
K\setminus F$ and assume that $t$ is not a $p$-th power in $K$. If $p(K,t)$ is
definable in $K$, then the set
\[%
\calB(K,t):=\{(t^{p^s}, x^{p^s}, x): s \in \Z_{> 0}, x \in K\}
\]%
is definable in $K$.

\end{corollary}%
\subsection{Constructing a model with
  addition and divisibility}
The final result we will need to construct a model of
$\langle\Z_{>0},+,\mid\rangle$ is the the following proposition.
\begin{proposition}%
  Let $t$ be as above. If $\calX(K)=\{ (x, x^{p^s}), x \in K, s \in
  \Z_{>0}\}$ is definable in $K$, then the set
\[
\calC(K,t):=\{(t^{p^{a}}, t^{p^{b}}, t^{p^{a+b}}): a, b > 0\}
\]%
is definable in $K$.
\end{proposition}%
\begin{proof}%
Consider the following system of equations:
\begin{equation}
\left \{
\begin{array}{c}
x-1=t^{p^{a}}\\%
\exists l\in \Z_{> 0}: z =((t+1)t^{p^{b}})^{p^l}\\%
\exists j\in \Z_{> 0}: z/x= t^{p^j}%
\end{array}%
\right .%
\end{equation}%
We claim that for any $a, b > 0$, if this system has solutions $x, z
\in K$ then $x/z=t^{p^{a+b}}$. Indeed, from the first equation we
conclude that $x=(t+1)^{p^a}$. From the second equation we get that
$z=(t+1)^{p^l}t^{p^{b+l}}$. Finally, from the third equation we have
that $(t+1)^{p^l-p^a}t^{p^{b+l}}=t^{p^j}$. The only way this equality
can hold is for $l=a$ and $j = a+b$. Conversely, we can always satisfy
the system if $x=(t+1)^{p^a}$ and $z=(t+1)^{p^a}t^{p^{b+a}}$.
\end{proof}%

We are now ready for the following theorem.
\begin{theorem}%
\label{thm:model}
Assume that $t$ has a simple pole or a simple zero.  Suppose that the set
$p(K,t)$ is definable in $K$. Then
$\langle\Z_{>0},+,\mid\rangle$ has a model over $K$.
\end{theorem}%
\begin{proof}%
  We map $s >0$ to $t^{p^s}$. Then $s=s_1+s_2 \Leftrightarrow
  (t^{p^{s_1}},t^{p^{s_2}},t^{p^s}) \in \calC(K,t)$. Further $s_1 \Big{|}
  s_2$ if and only if $(p^{s_1}-1)\Big{|} (p^{s_2}-1)$ if and only if
  there exists $x \in K$ such that
  \begin{equation}
  \label{eq:divide}
  x^{p^{s_1}-1}=t^{p^{s_2}-1},
  \end{equation}%
  since at least one pole or zero of $t$ is simple. Indeed suppose
  that the equality holds and let $\qq$ be a simple pole or zero of
  $t$. Then
\[%
p^{s_2}-1=\ord_{\qq}t^{p^{s_2}-1}= \ord_{\qq}x^{p^{s_1}-1} \equiv 0 \mod
p^{s_1}-1.
\]%
Conversely, if $p^{s_2}-1=l(p^{s_1}-1)$ for some $l \in \Z_{>0}$, then
we can set $x=t^l$ and (\ref{eq:divide}) will hold.

Hence $s_1
\Big{|} s_2$ if and only if

\[\exists \; x, y \in K \left((t^{p^{s_1}},y,x) \in
\calB(K,t) \;\land\; y/x=t^{p^{s_2}}/t\right).\]

 The result now follows from the fact that the sets
$p(K,t),\calB(K,t)$ and $\calC(K,t)$ are all definable in $K$.
\end{proof}%
\subsection{Defining $p$-th powers of one special element} We now
address the issue of defining {$p$-th} powers of one specific element
when the constant field is perfect. In the next proposition we observe
that if we avoid ramified zeros and poles and consider rational
functions only, we have the desired result.
\begin{proposition}%
\label{prop:version1}
Assume $t$ has no zeros or poles ramifying in the extension $K/F(t)$.
Let $r=1$ if $p >2$ and let $r=2$ if $p=2$. Assume that $F$ is perfect
and for some element $w \in F(t)$, having no poles or zeros at the
primes ramifying in the extension $K/F(t)$, there exists $u, v \in K$
such that the following system is satisfied.
\begin{equation}%
\label{eq:alreadydown}%
\left \{ \begin{array}{c}%
\frac{1}{t}-\frac{1}{w}=u^{p^r}-u\\
t-w=v^{p^r}-v%
\end{array} \right . %
\end{equation}%
Then for some $s \in \Z_{\geq 0}$ we have that $w=t^{p^{rs}}$.
Conversely, if $w=t^{p^{rs}}, s \geq 0$, then there exist $u, v \in F(t)$
satisfying (\ref{eq:alreadydown}). (For the last assertion we do not
need the requirement that $F$ is perfect.)
\end{proposition}%
\begin{proof}%
  Given our assumptions on $F$ and $w$ the proof of this proposition
  is identical to the proofs of Lemma 8.3.3, Corollary 8.3.4, and
  Proposition 8.3.8 of \cite{Sh34}.%
\end{proof}%
Unfortunately, we cannot always assume that an arbitrary rational
field element $w$ avoids all ramified poles and zeros. However, this
problem can be solved rather easily if we modify the
equations. The next remark and the proposition below deal with an
arbitrary element $w \in F(t)$.

\begin{remark}\label{prev}
 We recall from
Notation~\ref{not:K} that $K/F$ was a function field of positive
characteristic $p$ and $F_0$ denoted the algebraic closure of a
finite field in $F$.
By Proposition~\ref{le:convert} and Lemma~\ref{le:ramify} from the
appendix we can enlarge the constant field and assume that $F_0$
contains elements $c_0=0, c_1,\ldots,c_{2n(\alpha)+2}$
such that when $i\neq j$ we have for all $k \in
\Z_{\geq 0}$ that $c_i^{p^k} \not = c_j$. Here $n(\alpha)$ is the
constant that is defined in Lemma~\ref{le:ramify} from the appendix.
\end{remark}

We can now prove the following proposition.
\begin{proposition}%
\label{prop:below}%
Assume $F$ is perfect and $t$ is not a $p$-th power in $K$.
Let $c_0, \dots,c_{2n(\alpha)+2}$ be as in Remark~\ref{prev}, and 
let $V_i = \{c_i^{p^k}, k\in \Z_{\geq 0}\}$.

  Let $r=1$
if $p >2$ and let $r=2$ if $p=2$. Let $w \in F(t)$ and suppose that
for all $i\not = j \in \{1,\ldots, 2n(\alpha)+2\}$ for some $a \in
V_i, b\in V_j$ there exist $u_{i,j,a,b}, v_{i,j,a,b} \in K$ such that
\begin{equation}%
\label{eq:alreadydown1}%
\left \{ \begin{array}{c}%
\frac{t-c_i}{t-c_j}-\frac{w-a}{w-b}=u_{i,j,a,b}^{p^r}-u_{i,j,a,b}\\
\frac{t-c_j}{t-c_i} - \frac{w-b}{w-a}=v_{i,j,a,b}^{p^r}-v_{i,j,a,b}%
\end{array} \right . %
\end{equation}%
Then for some $s \in \Z_{\geq 0}$ we have that $w=t^{p^{rs}}$.
Conversely, if $w=t^{p^{rs}}$ for some $s \in \Z_{\geq 0}$ then the
equations can be satisfied as specified above (even if $F$ is not
perfect).
\end{proposition}%
\begin{proof}%
  First of all, using an argument similar to the one used in Lemma
  8.3.10 of \cite{Sh34}, we conclude that for some $c_i, c_j$ and $a
  \in V_i, b \in V_j$, we have that $t-c_i, t-c_j, w-a, w-b$ do not have
  zeros at any prime ramifying in the separable extension $K/F(t)$ and
  therefore $\displaystyle \frac{t-c_i}{t-c_j}, \frac{w-a}{w-b} \in
  F(t)$ do not have zeros or poles at any primes ramifying in the
  extension $K/F(t)$. Now applying Proposition \ref{prop:version1} we
  conclude that either $\displaystyle \frac{w-a}{w-b} = \left(
    \frac{t-c_i}{t-c_j}\right)^{p^{rs}}$, $s >0$ or $\displaystyle
  \frac{w-a}{w-b} = \frac{t-c_i}{t-c_j}$. In the first case we can
  take the $p^r$-th ``root'' of all our equations as in Lemma 8.3.1
  and Lemma 8.3.2 of \cite{Sh34}. In the second case we obtain $w =
  a_1t + a_2$ for some $a_1, a_2 \in F_0$. However, if we plug in this
  expression for $w$ into our equations with $c_0=0$ we obtain a
  contradiction unless $a_1=1$ and $a_2=0$. Thus, we conclude that
  $w=t^{p^{rs}}$. Finally, the satisfiability assertion follows as
  before from Proposition 8.3.8 of \cite{Sh34}.
\end{proof}%

To define $p$-th powers of a special element $t$ over fields of
transcendence degree one and higher transcendence degree we will use
some of the equations that were used in \cite{Sh15} and \cite{Eis}.
What we need to make the same arguments go through in our more general
setup is a set of equations over $K$ that forces its solutions to be
in the rational function field $F(t)$ and which are satisfied by all
elements $t^{p^s}$, $s \in \mathbb{Z}_{>0}$. I.e., we want a
set $\mathcal{S}$ which is definable in $K$ such that $p(K,t) \subseteq
\mathcal{S} \subseteq F(t)$ and thus we can apply Proposition
\ref{prop:below}. This will be accomplished in the next two sections.
For the transcendence degree one case, the equations defining
$\mathcal{S}$ are given in Corollary~\ref{finalpart1}. For higher
transcendence degree, they are given in Proposition~\ref{finalpart2}
below.

\section {Defining $p$-th powers for function
  fields whose constant field is algebraic}\label{td1}
 Let $K/F$ be a function field in one variable
  of positive characteristic $p$ with $F$ algebraic over a
  finite field. When $F$ has an extension of degree
  $p$, the results in
\cite{Sh15} and \cite{Eis} show that the existential theory of $K$ and
hence also the first-order theory of $K$ are undecidable.
Hence we can make additional assumptions
about the field of constants for the algebraic case and assume that we
are in a situation that is not covered by \cite{Sh15} or \cite{Eis}.
\begin{notation}%
\item Let $K/F$ be a function field in one variable
  of positive characteristic $p$.  Assume that $F$ is algebraic over a
  finite field and has no extension of degree $p$. 
\item Let $t$ be a
  fixed element of $K\setminus F$ which is not a $p$-th power in $K$.
\item We write $g_K$ for the genus of $K$, and when $f \in F[X,Y]$
  defines a plane curve $\mathcal{C}$ over $F$, we denote by $g_f$ the
  genus of the function field of $\mathcal{C}$ and also refer to
  this as the genus of $f$.
\end{notation}%
In this section we will show how to define $p$-th powers of the
element $t$ under the above assumptions.
\begin{lemma}%
\label{le:constsol}
For a pair of positive integers $k=p^l, u$, let
\[%
\displaystyle f_{k,u}(X,Y) = Y^{p^k}-Y +\frac{1}{\prod_{i=1}^u(X-c_i)}.
\]%
Then for any $a \in F$ there exists $b \in F$ such that $f_{k,u}(a,b)=0$.
\end{lemma}%
\begin{proof}%
  Fix an $a \in F$ and let $\alpha_1, \alpha_2$ be roots of
  $f_{k,u}(a,Y)$ in the algebraic closure of $F$. Then
  $(\alpha_1^{p^k} - \alpha_2^{p^k}) -(\alpha_1 -\alpha_2)=0$. Thus,
  $\alpha_1-\alpha_2=c$ is of degree $k=p^l$ over a field of $p$
  elements and therefore $c \in F$. Since $F$ is algebraic over a
  finite field, the extension $F(\alpha_1)/F$ is cyclic. Assume that
  $[F(\alpha_1):F]=m>1$ and let $\sigma \in \Gal(F(\alpha_1)/F)$ be a
  generator. Then for some $c \in F$ we have that $\sigma(\alpha_1) =
  \alpha_1 +c$ and $\id(\alpha_1) =\sigma^m(\alpha_1)=\alpha_1
  +mc=\alpha_1$. Thus $m\equiv 0\mod p$, and $F(\alpha)/F$ has a
  subextension of degree $p$ over $F$, contradicting our assumption on
  $F$.
\end{proof}%
\begin{lemma}%
\label{constantfield}%
There exists a set $A \subset K^2$, diophantine over $K$ such that $A \subset
F^2$ and for all $a \in F$ there exists $c \in F$ such that $(a,c) \in
A$.
\end{lemma}%
\begin{proof}%
  Before we proceed with the proof we should note that using the
  effective version Chebotarev Density Theorem (see \cite{Friednew},
  Proposition 6.4.8) one could show that any infinite field algebraic
  over a finite field is anti-Mordellic and therefore one could use a
  result of Poonen and Pop to see that $F$ is first-order definable in
  $K$ (see \cite{Popo}).  However in our case we can give a very
  simple existential definition of $F$ along the lines of
  \cite{Duret2}, \cite{Koenig} and \cite{Sh35}.

  The idea is to construct an equation $f$ whose genus is greater than
  the genus of $K$ and then use the Riemann-Hurwitz formula to show
  that all the $K$-rational solutions must be $F$-rational. We also
  have to ensure that $f$ has enough solutions over $F$.  Consider an
  equation $\displaystyle f_{k,u}(X,Y) = Y^{p^k}-Y
  +\prod_{i=1}^u\frac{1}{(X-c_i)}$, where $k$ and $u$ are as above and
  $c_1,\ldots, c_u$ are all distinct and in $F$. For sufficiently high $k$
  and $u$ the genus of this equation is higher than the genus of $K$.
  To see that this is so, assume $(u,p)=1$ and consider the field
  extension $F_{k,u}(X,Y)$ of $F(X)$ where $f_{k,u}(X,Y) =0$. It is
  clear that in this extension the primes corresponding to
  $(X-c_1),\ldots,(X-c_u)$ are completely ramified.  It is also clear from
  considering the difference between any two roots of this equation as
  in the lemma above, that no other prime of $F(X)$ is ramified in the
  extension $F_{k,u}(X,Y)/F(X)$.  Furthermore, the
  $F_{k,u}(X,Y)$-factor of $(X-c_i)$ is of relative degree 1 and also
  of degree 1 in $F(X,Y)$. Let $g_X=0$ be the genus of $F(X)$, and let
  $g_{f_{k,u}}$ be the genus of $f_{k,u}(X,Y)$. Then by the
  Riemann-Hurwitz formula and Remark 3.5.7 of \cite{Friednew}, we have
  that
\[%
2g_{f_{k,u}}-2 \geq p^k(g_X-2) + \deg \sum_{i=1}^{u}(p^k-1)\mte{P}_i,%
\]%
where for $i=1,\ldots,u$, we let $\mte{P}_i$ denote the prime above
$X-c_i$. Thus,%
\[%
g_{f_{k,u}} \geq \frac{1}{2}(u(p^k-1) - 2p^k+2)=\frac{1}{2}(p^ku
-u-2p^k+2)=\frac{(u-2)(p^k-1)}{2}.%
\]%
Now choose $k_0, u_0$ large enough so that $g_{f_{k_0,u_0}}$ is
greater than $g_K$, the genus of $K$. Let $f:=f_{k_0,u_0}$, let $F(X,Y)$
the corresponding field extension of $F(X)$, and
$g_f$ its genus.

Now assume that there exists a solution $x,y \in K\setminus F$ to
$f(X,Y)=0$. Then $F(X,Y) \isom F(x,y)$, so $F(X,Y)$ can be viewed as
a subfield of $K$. If $x$ is not a $p$-th power in $K$, then the
extension $K/F(x)$ is separable (see Chapter VI of \cite{Mason}) and
as a consequence of the Riemann-Hurwitz formula, $g_K \geq g_f$
contradicting the hypothesis.

If $x$ is a $p$-th power in $K$, then $\prod_{i=1}^u(x-c_i)$ is also a
$p$-th power in $K$ since $F$ is algebraic over a finite field, and
therefore all the coefficients $c_i$ of $f$ are also $p$-th powers.
Consequently $y$ is also a $p$-th power in $K$. Thus, by replacing all
the terms of $f$ by their $p$-th roots we can obtain a new equation
$f^{(1)}(X,Y)=0$ which is a ``$p$-th root'' of $f$.  The equation
$f^{(1)}$ has the same genus as $f$ because its genus only depends on
the values $k_0,u_0$ that were chosen above. Since $x$ and $y$ were
both $p$-th powers in $K$, the equation $f^{(1)}(X,Y)=0$ also has a
non-constant solution in $K$. Thus, at some point we will have an
equality $f^{(\ell)}(\tilde x, \tilde y)=0$, with the genus of
$f^{(\ell)}$ higher than the genus of $K$ and $\tilde{x}$ not a $p$-th
power in $K$. Consequently, $f(X,Y)=0$ can have constant solutions in
$K$ only. At the same time, by Lemma \ref{le:constsol}, for all $x \in
F$ we have $y \in F$ so that $f(x,y)=0$.
\end{proof}%

\begin{proposition}%
\label{prop:weak}
Suppose for some $w \in K$, for infinitely many primes $\mathfrak P$ of
$F(t)$ we have that
\begin{equation}%
\label{eq:equiv}
w\equiv a(\mathfrak P) \mod \mathfrak P,
\end{equation}%
where $a(\mathfrak P) \in F$. Then $w \in F(t)$.
\end{proposition}%
\begin{proof}%
  This proof follows from an argument similar to the argument in the
  proof of Theorem 10.1.1 of \cite{Sh34}. The main difference is that
  we do not assume that the prime $ \mathfrak P$ is inert in the
  extension $K/F(t)$. However, as long as the equivalence
  (\ref{eq:equiv}) holds for all the factors of the given prime below,
  the argument is unchanged.
\end{proof}%
\begin{corollary}%
\label{finalpart1}%
Suppose that for some $w \in K$ and infinitely many $(a, b) \in F^2$ we
have that the following system has a solution $u_a$ in $K$:
\begin{equation}
\label{eq:getdown}%
\frac{1}{t-a}-\frac{1}{w-b}=u_a^p-u_a
\end{equation}%
Then $w \in F(t)$.  Conversely, if for some positive integer $s$ we
have that $w=t^{p^s}$, then for any $a \in F$, there exist $b \in F, u_a
\in F(t)$ such that equation (\ref{eq:getdown}) is satisfied.
\end{corollary}%
\begin{proof}%
  First of all observe that the extension $K/F(t)$ is separable since
  $t$ is not a $p$-th power in $F(t)$. (See Lemma B1.32 of
  \cite{Sh34}.) Thus only finitely many primes ramify in the extension
  $K/F(t)$. Therefore for all but finitely many $a \in F$, for any
  factor $\pp_a$ of the rational prime $\mathfrak P_a$ which is the
  zero of $t-a$ in $F(t)$, it is the case that $\ord_{\pp_a}(t-a)=1$
  in $K$. On the other hand, for any pole $\qq$ of $u_a$ in $K$ we
  have that $\ord_{\qq_a}(u_a^p-u_a) \equiv 0 \mod p$. Thus, for all but
  finitely many $a \in F$, for all factors $\pp_a$ of the rational
  prime $\mathfrak P_a$ in $K$ we have that $\ord_{\pp_a}(w-b) >0$.
  In other words, for infinitely many $(a,b) \in F^2$ we have that $w
  \equiv b \mod \mathfrak P_a$, where $\mathfrak P_a$ is, as above, the
  zero divisor in $K$ and $F(t)$ of $t-a$. Now the first assertion of
  the corollary follows by Proposition \ref{prop:weak}. The second
  assertion of the corollary follows from Proposition 8.3.8 of
  \cite{Sh34}.
\end{proof}%

Finally note that we have assembled all the parts (i.e.\ Proposition
\ref{prop:below}, Lemma \ref{constantfield} and Corollary
\ref{finalpart1}) for the main result of this section.
\begin{theorem}%
  Let $K$ be a function field (in one variable) whose constant field
  $F$ is algebraic over a finite field of
  characteristic $p>0$.  Let $t$ be an element of $K\setminus F$ which
  is not a $p$-th power in $K$.  Then the set $p(K,t)=\{x \in K:
  \exists s \in \Z_{> 0} \, \, x =t^{p^s}\}$ is first-order definable
  in $K$.
\end{theorem}%

\section{Defining $p$-th powers over fields of higher transcendence
  degree}\label{highertd}
Let $K$ be a function field of characteristic $p>2$ with constant
field $F$, and assume that $F$ has transcendence degree at least one
over a finite field.  To define $p$-th powers of a suitable element
$t$ we will use a theorem by Moret-Bailly (\cite{MB3}). Here, we
quickly review his notation and state the theorem in the form we need.
\begin{definition}
Let $u:A \to B$ be a morphism of abelian groups. We say that $u$ is
{\em almost bijective} if $u$ is injective and Coker $u$ is a finite $p$-group.
\end{definition}
By \cite[Theorem 1.8]{MB3}, the following theorem holds:
\begin{theorem}
\label{thm:MB}
Let $F$ be a field of characteristic $p>2$, and assume that $F$
contains an element which is transcendental over $\mathbb{F}_p$.
Let $K$ be a function field in one variable with constant field $F$, and
let $$E:y^2=P(x)$$ be an elliptic curve which is defined over a finite
field contained in $F$.  There exists a non-constant element $t \in K$
such that $t$ is not a $p$-th power in $K$ and the elliptic curve
$\calE$ given by
  $$\calE: P(t)\,y^2 = P(x)$$ has the property that the natural
  homomorphism $\calE(F(t))\hookrightarrow \calE(K)$ induced by the
  inclusion $F(t) \hookrightarrow K$ is almost bijective.
\end{theorem}

\begin{notation}%
\label{notation}%
From now on, let $P(x),E$ and $t$ be as in Theorem~\ref{thm:MB}. Let
$s$ be an element in a quadratic extension of $K$ satisfying $s^2=
P(t)$.  Let $q=p^r$ be the size of a finite field containing all the
coefficients of the equation defining $E$.

Let $F'$ be an algebraic closure of $F$, and let $K'=F'K$.
\end{notation}
By \cite[Theorem 1.8]{MB3}, it follows that the natural homomorphism
$\calE(F(t)) \hookrightarrow \calE(K')$ is still almost bijective.
\begin{proposition}
The set $\calE(F(T))$ is diophantine over $K$ and over $K'$.
\end{proposition}
\begin{proof}
  Let $A:= \calE(F(T))$ and $B:= \calE(K)$. The set $B$ is
  clearly diophantine over $K$. By Theorem~\ref{thm:MB}, $A$ is a
  subgroup of finite index in $B$ and $B/A$ is a finite $p$-group.

  Hence for some integer $k$ we have that $p^k B \subseteq A$ and $p^kB$ has
  finite index in $B$. Since $B$ is diophantine over $K$, and since
  multiplication by $p^k$ is given by explicit equations, the set $p^k
  B$ is diophantine over $K$. It is easy to see that this implies that
  $A$ is diophantine over $K$:

Let $Q_1,\dots,Q_{\ell}$ be coset representatives for $p^kB$ in
$A$. Then for $P \in \calE(K)$
\[
P \in A \Leftrightarrow (\exists\, S \in
p^kB)(P=S+Q_1) \lor \dots \lor (P=S+Q_{\ell}).
\]
The same argument with $K$ replaced by $K'$ shows that $A$ is also
diophantine over $K'$.
\end{proof}
From the proposition above we also obtain the following easy corollary.
\begin{corollary}
  There exists a polynomial equation $R(u,v,x_1,\ldots,x_l) \in
  K[u,v,x_1,\ldots,x_l]$ such that $R(u,v,x_1,\ldots,x_l)=0$ for some
  $u,v,x_1,\ldots,x_l \in K'$ implies $(u,v)$ are affine coordinates of a
  point in $\calE(F(t))$.  Conversely, if $(u,v)$ are affine
  coordinates of a point in $\calE(F(t))$ the equation
  $R(u,v,x_1,\ldots,x_l)=0$ can be satisfied with $x_1,\ldots, x_l \in K$.
\end{corollary}

Next we observe that $p$-th powers occur as affine coordinates of
points of $\calE$.
\begin{lemma}%
\label{le:powersoft2}
  Let $\calE, s, t, q$ be as in Notation~\ref{notation}. The point
  $(t^{q^m}, s^{q^m-1}) \in \calE(F(t))$.
\end{lemma}%
\begin{proof}%
  Observe that $(P(t))^{q^m}= P(t^{q^m})$.  Thus, $P(t)(s^{q^m-1})^2=
  (P(t))^{q^m}=P(t^{q^m})$.  Also $q^m -1$ is
  even, so the point $(t^{q^m}, s^{q^m-1})$ has coordinates in the ground
  field.
\end{proof}%
We conclude with the proposition defining $p$-th
powers of $t$ for the case of $K$ of transcendence
degree greater than one.
\begin{proposition}\label{finalpart2}%
  Assume that for some $z, w, u, v \in K$ the following system is
  satisfied over $K'$.%
 \begin{equation}
\label{eq:adec}%
\left \{ \begin{array}{c}%
R(w,z, x_1,\ldots,x_l)=0\\
\forall i, j \in \{1,\ldots,2n(\alpha)+2\} \exists a \in V_i, b\in V_j:\\
\frac{t-c_i}{t-c_j}-\frac{w-a}{w-b}=u_{i,j,a, b}^p-u_{i,j,a, b}\\
\frac{t-c_j}{t-c_i} - \frac{w-b}{w-a}=v_{i,j,a,b}^p-v_{i,j,a,b}%
\end{array} \right . %
\end{equation}%
Then for some $s \in \Z_{\geq 0}$ we have that $w=t^{p^s}$. Conversely,
if $w=t^{p^s}$, then the system has solutions in $K$.
\end{proposition}%
Finally we note that $t$ selected so that Theorem \ref{thm:MB} holds
might have zeros and poles which are not simple.  Let $t' =
\frac{t-a}{t-b}$ be such that that all of its poles and zeros in $K$
are simple.  Observe that $F(t) = F\left(\frac{t-a}{t-b}\right)$ and we can
generate $p$-th powers of $t'$.  Thus in Theorem \ref{thm:model}
and Proposition \ref{prop:needt} we can replace $t$ by $t'$ if necessary.%
\begin{remark}
  The proof of Proposition~\ref{prop:needt} that $p$-th powers of the
  special element $t$ allow us to define $p$-th powers of arbitrary
  elements $x$ in $K$ only used equations involving existential
  quantifiers. The proof of Proposition~\ref{finalpart2} which defined
  $p$-th powers of $t$ also only used existential quantifiers. So when
  the constant field of $K$ contains transcendental elements over
  $\F_p$, the set $\{(x,x^{p^s}):x \in K, s \in \Z_{\geq 0}\}$ is
  actually {\em existentially} definable in $K$ .
\end{remark}

\section{Ring language}\label{coeff}

In this section we address the issue of the language needed to produce
an undecidable set of sentences.  We have already shown that we can
construct a model of the positive integers $\Z_{>0}$ (and hence also
of $\Z_{\geq 0}$) in the function field $K$.  We used this model to
prove Theorem~\ref{maintheorem}.

In this section we will show that this easily implies
Theorem~\ref{maintheorem2} by using a result of R.\ Robinson, so we
obtain undecidability of the first-order theory of $K$ in the language
of rings {\em without} parameters.

\begin{corollary}[Theorem~\ref{maintheorem2}]\label{noparameters}
  Let $K$ be a function field of characteristic $p> 2$. Then the
  first-order theory is undecidable in the language of rings without
  parameters. When $K$ is a function field in one variable whose
  constant field is algebraic over a finite field, then we also obtain
  undecidability in characteristic 2.
\end{corollary}
\begin{proof}
  From the previous sections it follows that the equations we used to
  construct a model of $\Z_{\geq 0}$ are in the language $L_{\bar d}=
  \langle +, \cdot\,;\; 0,1, \{ d_1, \dots, d_r \} \rangle$, for fixed
  elements $d_1, \dots, d_r$ of $K$. In other words, we are working in
  the language of rings with finitely many parameters. To show that we
  can achieve undecidability in the ring language {\it without
    parameters}, we use a result of R.\ Robinson who gave an example
  of a finitely axiomatizable and essentially undecidable theory $Q$
  (\cite[p.\ 32]{TMR53}).  A theory is {\em essentially undecidable}
  if any consistent extension of it is also undecidable.  Since $Q$ is
  a subtheory of $\Z_{\geq 0}$ (\cite[p.\ 51]{TMR53}), the axioms of
  $Q$ hold in $\Z_{\geq 0}$. Let $Ax(\Z_{\geq 0})$ be the conjunction
  of all the axioms of $Q$. For a sentence $\psi$ in the language $L=
  \langle +, \cdot, 0,1 \rangle$ let $\phi_K(\psi,\bar d)$ be its
  translation in our model, and consider the set of all $L$-sentences
  $\psi$ for which
\begin{equation}
\label{forall}%
\forall \bar w (\phi_K(Ax(\Z_{\geq 0}),\bar w) \to \phi_K(\psi,\bar w))%
\end{equation}%
is true in $K$. This set contains the axioms of $Q$ and therefore the
theory generated by these sentences is an extension of $Q$. The
extension is consistent.  Suppose not. Then for some $\psi$ as above we
have that $Ax(\Z_{\geq 0}) \to \lnot \psi$ holds in $\Z_{\geq 0}$, and hence
$\phi_K(Ax(\Z_{\geq 0}),\bar d) \to \lnot \phi_K(\psi,\bar d)$ holds in $K$.  But
this contradicts (\ref{forall}).  Since the collection of all
$L$-formulas $\psi$ satisfying (\ref{forall}) in $K$ is undecidable, the
set of all formulas of the form (\ref{forall}) true in $K$ is also
undecidable.  Finally we note that the formulas in this set are in
$L$.
\end{proof}
\section{Open questions}
Even though we proved that the first-order theory of function fields
of characteristic $p>2$ is undecidable, we needed transcendental
elements in the construction of the model of the nonnegative integers.
So the following questions arise naturally:
\begin{question}
  Let $K$ be a function field as in Theorem~\ref{maintheorem}. Does
  $K$ admit a model of $\langle \Z_{\geq 0},+,\cdot \rangle$ in which
  the equations defining the model have integer coefficients?
\end{question}
\begin{question}
Is the degree of unsolvability of the first-order theory of $K$ at
least that of the first-order theory of $\Z$? 
\end{question}

\section{Appendix}\label{appendix}
 The following proposition and lemma are used to handle
 constant field extensions of function fields.
\begin{proposition}%
\label{le:convert}%
Let $Q_i$ be either ``$\forall$'' or ``$\exists$''.  Let $M/K$ be a
finite extension of fields with $M$ not algebraically closed.    Let
$P(t_1,\ldots,t_r,x_1,\ldots,x_k) \in M_0[t_1,\ldots,t_r,x_1,\ldots,x_k]$.  Let
\begin{equation}%
\label{sen:1}%
A_M=\{(t_1,\ldots,t_r) \in M^r: Q_1 x_1\in M \ldots Q_k x_k \in M: P(t_1,\ldots, t_r,x_1,\ldots,x_k)=0 \}
\end{equation}%
be a first-order definable set.
Then there exists a polynomial $T(u_1,\ldots,u_m, y_1,\ldots,y_l) \in K_0[u_1,\ldots,u_m,
y_1,\ldots,y_l]$,  and a first-order definable set
\begin{equation}%
\label{sen:2}%
A_K=\{(u_1,\ldots, u_m) \in K^m: Q_{k+1} y_1 \in K\ldots Q_{k+l}  y_l \in K :
T(u_1,\ldots, u_m,y_1,\ldots,y_l)=0
\end{equation}%
such that for any $(t_1,\ldots, t_r) \in M^r$ we have that $(t_1,\ldots, t_r)
\in A_M$ if and only if for some $m$-tuple $(u_1,\ldots,u_m) \in K^m$ we
have that $(u_1,\ldots,u_m) \in A_K$.  Thus, if $M$ has a first-order model
of $\Z$ in the language of rings augmented by finitely many parameters
from $M$, then $K$ has a first-order model of $\Z$ in the language of
rings with finitely many parameters from $K$.
\end{proposition}%
The proof of the proposition requires standard ``rewriting''
techniques utilizing a basis of $M$ over $K$ and the fact that over a
field which is not algebraically closed we can replace a finite set of
equations by a single equivalent equation.

Proposition \ref{le:convert} will play a role in case we need to
extend the field of constants to ensure that we have ``enough''
conjugacy classes of constants algebraic over a finite field relative
to the number of primes ramifying in $K/F(t)$ or $K'/F'(t)$, where
$F'$ is the algebraic closure of $F$ and $K'=F'K$.  In this connection
we have the following lemma.

\begin{lemma}%
\label{le:ramify}%
Let $\alpha$ be any generator of $K$ over $F(t)$ with $K/F(t)$ separable.
Let $h(T)=a_0 + a_1T+ \ldots T^{n}$ be the monic irreducible polynomial of
$\alpha$ over $F(t)$.  Let $D(\alpha)= {\mathbf N}_{K/F(t)}(h'(\alpha))$ where
$h'(T)$ is the derivative of $h(T)$ with respect to $T$.  Let $P(\alpha)$
be the pole divisor of $\prod_{i=0}^{n-1}a_i$.  Since $F(t)$ is a
rational function field, $D(\alpha)$ and $P(\alpha)$ are both polynomials in
$t$.  Let $n(\alpha)$ be the degree of the polynomial $D(\alpha)P(\alpha)$.  Let
$\hat{F}$ be any algebraic extension of $F$.  Then the number of
$\hat{F}(t)$ primes ramifying in the extension $\hat{F}K/\hat{F}(t)$
is less or equal to $n(\alpha)$.
\end{lemma}%
\begin{proof}%
  Since there is no constant field extension in the extension $K/F(t)$
  (see Notation \ref{not:K}), we have that $K$ and $\hat{F}(t)$ are
  linearly disjoint over $F(t)$.  Thus, $\hat{F}K=\hat{F}(t)(\alpha)$ and
  $[\hat{F}K:\hat F(t)]=[K:F(t)]=n$.  Therefore if $\qq$ is a prime of
  $\hat F(t)$ ramified in the extension $\hat FK/\hat F(t)$ we have
  two options: either $\qq$ divides the discriminant $D(\alpha)$ of the
  power basis of $\alpha$ or $\alpha$ is not integral at $\qq$ and $\qq$
  divides $P(\alpha)$.  In either case $\qq$ divides $D(\alpha)P(\alpha)$. Since
  $P(\alpha)D(\alpha)$ is a polynomial, its degree is invariant under any
  constant field extension and therefore the number of primes dividing
  $P(\alpha)D(\alpha)$ in $\hat F$ is bounded by the degree of this
  polynomial.
\end{proof}%

\section*{Acknowledgments}
The authors would like to thank Arno Fehm and Stephen Simpson for
helpful discussions leading to Corollary~\ref{noparameters}.

\end{document}